\documentclass[12pt,oneside,english]{amsart}
\usepackage[latin1]{inputenc}
\usepackage{amssymb}

\makeatletter

\newtheorem{remark}{Remark}

\usepackage{babel}

\makeatother
\begin{document}
\baselineskip=17pt
\title[On a Luschny question]{On a Luschny question}

\author{Vladimir Shevelev}
\address{ Department of Mathematics \\Ben-Gurion University of the
 Negev\\Beer-Sheva 84105, Israel. e-mail:shevelev@bgu.ac.il}


\begin{abstract}
 Let $E_n(x)$ be Euler polynomial, $\nu_2(n)$ be $2-$adic order of $n,$ $\{g(n)\}$ be the characteristic sequence for $\{2^{n}-1\}_{n\geq1}.$ Recently Peter Luschny asked (cf. \cite{5}, sequence A135517): is A135517(n) the denominator of $E_n(x) - E_n(1)?$  According to a formula in A091090, this question is equivalent to the following one: is the denominator of $E_n(x) - E_n(1)$ equal to $2^{\nu_2(n+1)-g(n)}?$ In this note we answer this question in the affirmative. 
\end{abstract}

\maketitle

\section{Introduction}
Let $E_n(x)$ be Euler polynomial, $\nu_2(n)$ be $2-$adic order of $n,$ $\{g(n)\}$ be the characteristic sequence for $\{2^{n}-1\}_{n\geq1}.$ Recently  Peter Luschny asked (cf. \cite{5}, sequence A135517): is A135517(n) the denominator of $E_n(x) - E_n(1)?$  According to a formula in A091090, this question is equivalent to the following one: is the denominator of $E_n(x) - E_n(1)$ equal to $2^{\nu_2(n+1)-g(n)}?$ In this note we answer this question in the affirmative. Our proof is based on finding a simple explicit expression for the coefficients of Euler polynomial.
\begin{remark}\label{r1}
 Note that Peter Luschny published in OEIS sequence $A290646$ in which he for the first time asked his question: $"Is\enskip A290646 = A135517?"$ and in the version 2 of this note we referred to $A290646.$ But after publication of version 2, the sequence $A290646$ became a replacement sequence with no relation to this topic. Therefore, in order not to cause any inconvenience to the readers, we are forced to give this third version.
\end{remark}

\section{Several classic formulas and theorems}
Euler polynomials $E_n(x)$ are defined by generating function
\begin{equation}\label{1}
\frac{2e^{xt}}{e^t+1} = \sum_{n=0}^{\infty} E_n(x)\frac{t^n}{n!}.
\end{equation}
Below we use several known relations \cite{1}
\begin{equation}\label{2}
(-1)^nE_n(-x)=2x^n-E_n(x);
\end{equation}
\begin{equation}\label{3}
E_n(0)=-E_n(1)=-\frac{2}{n+1}(2^{n+1}-1)B_{n+1}, \enskip n=1,2,...,
\end{equation}
where $\{B_n\}$ are Bernoulli numbers;
\begin{equation}\label{4}
E'_n(x)=nE_{n-1}(x).
\end{equation}
We use also the formula which is obtained by combining formulas (14) and (18) in \cite{2} (see also \cite{6}):
\begin{equation}\label{5}
B_n=\frac{n}{2(2^n-1)}\sum_{j=0}^{n-1}(-1)^jS(n,j+1)\frac{j!}{2^j},
\end{equation}
where $\{S(n,j)\}$ are the Stirling numbers of the second kind.\newline
\indent Further recall that, according to Von Staudt-Clausen theorem \cite{4},
we have
\begin{equation}\label{6}
B_{2n}=I_n - \sum \frac{1}{p},
\end{equation}
where $I_n$ is an integer, $\{p\}$ are primes for which $p-1$ divides $2n.$\newline
Finally, denote by $t(n,k)$ the number of carries which appear in addition
$k$ and $n-k$ in base 2, or, the same, in subtracting $k$ from $n.$ Then, by
 Kummer's known theorem (cf.\cite{3}),
\begin{equation}\label{7}
2^{t(n,k)}||\binom{n}{k},
\end{equation}
i.e., $t(n,k)$ is 2-adic order of $\binom{n}{k}.$

\section{Explicit formula for coefficients of $E_n(x)$}
Let
$$E_n(x)=e_0(n)x^n+e_1(n)x^{n-1}+e_2(n)x^{n-2}+...+e_{n-1}(n)x+e_n(n). $$
Using (\ref{2}), we immediately find
$$ e_0(n)=1, e_2(n)=e_4(n)=...=0.$$
So, we have
\begin{equation}\label{8}
E_n(x)= x^n + \sum_{odd\enskip k=1,...,n} e_k(n)x^{n-k}.
\end{equation}

Further, by (\ref{4}), $e_k(n)$ satisfies the difference equation
\begin{equation}\label{9}
e_k(n)=\frac{n}{n-k}e_k(n-1).
\end{equation}
It is easy to see that the solution of (\ref{9}) is
\begin{equation}\label{10}
e_k(n)=C_k\binom{n}{k}.
\end{equation}
Firstly, let us find $e_k(n)$ for odd $n.$ Then by (\ref{3})
 $$e_n(n)=C_n=E_n(0)=-\frac{2}{n+1}(2^{n+1}-1)B_{n+1}.$$
So, by (\ref{10}) for odd $n$ we find
\newpage
\begin{equation}\label{11}
e_k(n)=-\frac{2}{k+1}(2^{k+1}-1)B_{k+1}\binom{n}{k}.
\end{equation}
Now let $n$ be even. Let us show that the formula for $e_k(n)$ does not change.\newline
 \indent Indeed, again by (\ref{4}), we have
$$(n+1)E_n(x)=E'_{n+1}(x).$$
So, by (\ref{9}) and (\ref{11}), we have
$$x^n + \sum_{odd\enskip k=1,...,n} e_k(n)x^{n-k}=$$
$$\frac{1}{n+1}((n+1)x^n + \sum_{odd\enskip k=1,...,n+1}C_k\binom{n+1}{k}(n+1-k)x^{n-k})=$$
$$x^n+\sum_{odd\enskip k=1,...,n-1}C_k\binom{n+1}{k}\frac{n+1-k}{n+1}x^{n-k}. $$
Hence, for even $n$ we find
\begin{equation}\label{12}
e_k(n)=-\frac{2}{k+1}(2^{k+1}-1)B_{k+1}\binom{n}{k},
\end{equation}
that coincides with (\ref{11}).\newline
\indent  Let $x$ be a rational number. Below we denote by $N(x)$ the numerator of $x$ and by
$D(x)$ the denominator of $x,$ such that $N(x)$ and $D(x)$ are relatively prime.\newline
Now note, that by (\ref{7}), $D(B_{k+1})$ ($k$ is odd) is an even square-free
number, while $N(B_{k+1})$  is odd. Hence, $(2^{k+1}-1)N(2B_{k+1})$ is odd number.
Finally, by (\ref{11})-(\ref{12}) and (\ref{5}) we have (for odd k):
\begin{equation}\label{13}
e_k(n)=-\binom{n}{k}\sum_{j=0}^k (-1)^jS(k+1, j+1)\frac{j!}{2^j}.
\end{equation}
This yields that $D(e_k(n))$ could be only a power of 2. This means
for (\ref{11})-(\ref{12}), that $D(e_k(n))$ is really
\begin{equation}\label{14}
D(e_k(n))= 2^{\nu_2(k+1)-\nu_2(\binom{n}{k})},
\end{equation}
where $\nu_2(n)$ is $2-$adic order of $n.$\newline
\indent Add that, since $sign(B_{2n})=(-1)^{n-1},$ then, by (\ref{11})-(\ref{12}), $sign(e_k(n))=(-1)^{\frac{k+1}{2}}.$

\section{Answer in the affirmative on the Peter Luschny question}
Let, according to the question,
$$E^*_n(x)=E_n(x)-E_n(1).$$
\newpage
By (\ref{3}),
\begin{equation}\label{15}
E^*_n(x)=E_n(x)+E_n(0).
\end{equation}
In case of odd $n,$ when $E_n(0)\neq0,$ $E^*_n(0)=2E_n(0).$ So, the formula
(\ref{14}) for the corresponding coefficients $e^*_k(n)$ takes the form
\begin{equation}\label{16}
D(e^*_k(n))= 2^{\nu_2(k+1)-\nu_2(\binom{n}{k})- \delta_{n,k}},
\end{equation}
where $\delta_{n,k}=1,$ if $k=n,$ and $\delta_{n,k}=0,$ otherwise.
Further, by (\ref{16}), we have

$$D(e^*_k(n))=
 2^{\nu_2((k+1)/\binom{n}{k})- \delta_{n,k}}=
 2^{\nu_2((n+1)/\binom{n+1}{k+1})- \delta_{n,k}}=$$
\begin{equation}\label{17}
2^{\nu_2(n+1)-\nu_2(\binom{n+1}{k+1})- \delta_{n,k}}.
\end{equation}

Now, by (\ref{7}), we have

$$D(E^*_n(x))= 2^{\max(odd \enskip k=1,...,n)(\nu_2(n+1)-t(n+1,k+1)- \delta_{n,k})}=$$
\begin{equation}\label{18}
2^{\nu_2(n+1)-\min(odd \enskip k=1,...,n)(t(n+1,k+1)+ \delta_{n,k})}.
\end{equation}

Let firstly $n=2^m-1,\enskip m\geq1.$ Then in (\ref{18}) we obtain the minimum in case $k=n$ when $t(n+1,k+1)=0.$ So, since $delta(n,n)=1,$ the minimum is $1=g(n).$
So, by (\ref{18}),
$$D(E^*_n(x))=2^{\nu_2(n+1)-g(n)}.$$
Let now a positive $n$ have not a form $2^m-1.$ Let us show that
in this case the minimum in (\ref{18}) is 0. Indeed, take $k=2^{u(n)-1}-1,$ where $u(n)$ is
the number of $(0,1)-$digits in the binary expansion of $n.$ Then $k<n$ and evidently $t(n+1,k+1)=0.$ Since
also $\delta(n,k)=0,$ then the minimum in (\ref{18}) is 0. So,
$$D(E^*_n(x))=2^{\nu_2(n+1)}$$
and since here $g(n)=0,$ we complete the proof.
\newline\newline

 \bfseries Acknowledgement. \mdseries The author thanks Jean-Paul Allouche for sending article \cite{2} and for useful discussions.

\end{document}